\documentclass{article}%
\usepackage{graphicx}
\usepackage[intlimits]{amsmath}
\usepackage{latexsym}
\usepackage{amsfonts}
\usepackage{amssymb}%
\makeatletter
\newfont{\footsc}{cmcsc10 at 8truept}
\newfont{\footbf}{cmbx10 at 8truept}
\newfont{\footrm}{cmr10 at 10truept}
\renewcommand{\ps@plain}{\renewcommand{\@oddfoot}{\footsc the electronic journal of combinatorics
{\footbf XX} (200X), \#R00\hfil\footrm\thepage}} \makeatother
\newtheorem{theorem}{Theorem}

\newtheorem{lemma}[theorem]{Lemma}

\newtheorem{problem}[theorem]{Problem}

\newenvironment{proof}[1][Proof]{\noindent{\textbf {#1}  }}  {\hfill$\Box$\bigskip}

\begin{document}

\title{The smallest eigenvalue of $K_{r}$-free graphs}
\author{Vladimir Nikiforov\\Department of Mathematical Sciences, University of Memphis, Memphis TN 38152, USA}
\date{{\small Submitted: 2004; Accepted: 200X; Published: 200X}\\
{\small MR Subject Classifications: 05C35}}
\maketitle

\begin{abstract}
Let $G$ be a $K_{r+1}$-free graph with $n$ vertices and $m$ edges, and let
$\mu_{n}\left(  G\right)  $ be the smallest eigenvalue of its adjacency
matrix. We show that
\[
\mu_{n}\left(  G\right)  <-\frac{2^{r+1}m^{r}}{rn^{2r-1}}.
\]

\end{abstract}

Our notation and terminology are standard (see, e.g. \cite{Bol98}). If $G$ is
a graph of order $n,$ we write $e\left(  X\right)  $ for the number of edges
induced by a set $X\subset V\left(  G\right)  ,$ $\Gamma\left(  u\right)  $
for the set of neighbors of a vertex $u\in V\left(  G\right)  ,$ and $\mu
_{n}\left(  G\right)  $ for the smallest eigenvalue of the adjacency matrix of
$G$.

Bollob\'{a}s and Nikiforov \cite{BoNi04} observed that if $G$ is a dense
$K_{r}$-free graph then $\mu_{n}\left(  G\right)  <-\alpha n$ for some
$\alpha>0,$ independent of $n$. In this note we prove the following more
precise statement.

\begin{theorem}
\label{th2} Suppose $r\geq2$ and $G$ is a $K_{r+1}$-free graph with $n$
vertices and $m$ edges. Then
\[
\mu_{n}\left(  G\right)  <-\frac{2^{r+1}m^{r}}{rn^{2r-1}}.
\]

\end{theorem}

To prove this theorem we need two preliminary results. Our main tool will be
the following inequality proved in \cite{BoNi04}: If $V\left(  G\right)
=V_{1}\cup V_{2}$ is a bipartition of the vertices of a graph $G$ then%
\begin{equation}
\mu_{n}\left(  G\right)  \leq\frac{2e\left(  V_{1}\right)  }{\left\vert
V_{1}\right\vert }+\frac{2e\left(  V_{2}\right)  }{\left\vert V_{2}\right\vert
}-\frac{2e\left(  G\right)  }{\left\vert V\left(  G\right)  \right\vert }.
\label{mainin}%
\end{equation}

Write $t\left(  G\right)  $ for the number of triangles of a graph $G$ and
$t^{\prime\prime}\left(  G\right)  $ for the number of its induced subgraphs
of order $3$ and size $1.$ For every vertex $u,$ let $t\left(  u\right)
=e\left(  \Gamma\left(  u\right)  \right)  $ and $t^{\prime\prime}\left(
u\right)  =e\left(  V\left(  G\right)  \backslash\Gamma\left(  u\right)
\right)  .$ Observe the following simple equalities
\[
3t\left(  G\right)  =\sum_{u\in V\left(  G\right)  }t\left(  u\right)  \text{
\ \ and \ \ }t^{\prime\prime}\left(  G\right)  =\sum_{u\in V\left(  G\right)
}t^{\prime\prime}\left(  u\right)  .
\]

Inequality (\ref{mainin}) implies the following lemma.

\begin{lemma}
\label{le3}Let $G$ be a graph with $n$ vertices and $m$ edges with no isolated
vertices. Then
\begin{equation}
\mu_{n}\left(  G\right)  \leq\frac{2n}{\left(  n^{2}-2m\right)  }\sum_{u\in
V\left(  G\right)  }\frac{t\left(  u\right)  }{d\left(  u\right)  }%
-\frac{4m^{2}}{n\left(  n^{2}-2m\right)  }. \label{ineq1}%
\end{equation}

\end{lemma}

\begin{proof}
We start by recalling the equality%
\begin{equation}
3t\left(  G\right)  =\sum_{u\in V\left(  G\right)  }d^{2}\left(  u\right)
-nm+t^{\prime\prime}\left(  G\right)  \label{maineq}%
\end{equation}
whose proof we shall outline for convenience. For every edge $uv\in E\left(
G\right)  $ we have
\[
\left\vert \Gamma\left(  u\right)  \cap\Gamma\left(  v\right)  \right\vert
=\left\vert \Gamma\left(  u\right)  \right\vert +\left\vert \Gamma\left(
v\right)  \right\vert -\left\vert \Gamma\left(  u\right)  \cup\Gamma\left(
v\right)  \right\vert =d\left(  u\right)  +d\left(  v\right)  -n+\left\vert
\overline{\Gamma\left(  u\right)  }\cap\overline{\Gamma\left(  v\right)
}\right\vert .
\]
Summing this equality over all $uv\in E\left(  G\right)  $ we obtain%
\[
3t\left(  G\right)  =\sum_{uv\in E\left(  G\right)  }\left(  d\left(
u\right)  +d\left(  v\right)  -n+\left\vert \overline{\Gamma\left(  u\right)
}\cap\overline{\Gamma\left(  v\right)  }\right\vert \right)  =\sum_{u\in
V\left(  G\right)  }d^{2}\left(  u\right)  -nm+t^{\prime\prime}\left(
G\right)  ,
\]
as claimed.

For every $u\in V\left(  G\right)  $ and bipartition $V_{1}=\Gamma\left(
u\right)  $, $V_{2}=V\left(  G\right)  \backslash\Gamma\left(  u\right)  ,$
inequality (\ref{mainin}) implies
\[
\mu_{n}\left(  G\right)  \leq\frac{2e\left(  V_{1}\right)  }{d\left(
u\right)  }+\frac{2e\left(  V_{2}\right)  }{n-d\left(  u\right)  }-\frac
{2m}{n}=\frac{2t\left(  u\right)  }{d\left(  u\right)  }+\frac{2t^{\prime
\prime}\left(  u\right)  }{n-d\left(  u\right)  }-\frac{2m}{n},
\]
and therefore,
\[
\mu_{n}\left(  G\right)  \left(  n-d\left(  u\right)  \right)  \leq
\frac{2t\left(  u\right)  }{d\left(  u\right)  }\left(  n-d\left(  u\right)
\right)  +2t^{\prime\prime}\left(  u\right)  -\frac{2m}{n}\left(  n-d\left(
u\right)  \right)  .
\]
Summing this inequality for all $u\in V\left(  G\right)  $, in view of
(\ref{maineq}), we obtain
\begin{align*}
\mu_{n}\left(  G\right)  \left(  n^{2}-2m\right)   &  \leq n\sum_{u\in
V\left(  G\right)  }\frac{2t\left(  u\right)  }{d\left(  u\right)  }-6t\left(
G\right)  +2t^{\prime\prime}\left(  G\right)  -\frac{2m}{n}\left(
n^{2}-2m\right) \\
&  =n\sum_{u\in V\left(  G\right)  }\frac{2t\left(  u\right)  }{d\left(
u\right)  }+2nm-2\sum_{u\in V\left(  G\right)  }d^{2}\left(  i\right)
-\frac{2m}{n}\left(  n^{2}-2m\right) \\
&  \leq2n\sum_{u\in V\left(  G\right)  }\frac{t\left(  u\right)  }{d\left(
u\right)  }-\sum_{u\in V\left(  G\right)  }d^{2}\left(  u\right)  \leq
2n\sum_{u\in V\left(  G\right)  }\frac{t\left(  u\right)  }{d\left(  u\right)
}-\frac{4m^{2}}{n},
\end{align*}
completing the proof.
\end{proof}

\begin{proof}
[Proof of Theorem \ref{th2}]Assume that $G$ has no isolated vertices; the
general case follows immediately. Our proof is by induction on $r.$ If $G$ is
triangle-free, Lemma \ref{le3} implies%
\[
\mu_{n}\left(  G\right)  \leq\frac{-4m^{2}}{n\left(  n^{2}-2m\right)  }%
<-\frac{4m^{2}}{n^{3}},
\]
so the assertion holds for $r=2;$ assume it holds for $r-1\geq2.$ Since, for
every $u\in V\left(  G\right)  ,$ the graph induced by $\Gamma\left(
u\right)  $ is $K_{r}$-free, the induction hypothesis implies%
\[
\mu_{n}\left(  G\right)  <-\frac{2^{r}t^{r-1}\left(  u\right)  }{\left(
r-1\right)  d^{2r-3}\left(  u\right)  }.
\]
Summing this inequality for all $u\in V\left(  G\right)  ,$ we obtain
\begin{equation}
n\mu_{n}\left(  G\right)  <-\frac{2^{r}}{r-1}\sum_{u\in V\left(  G\right)
}\frac{t^{r-1}\left(  u\right)  }{d^{2r-3}\left(  u\right)  }. \label{ineq2}%
\end{equation}
By H\"{o}lder's inequality, we find that%
\[
\left(  \sum_{u\in V\left(  G\right)  }d\left(  u\right)  \right)
^{\frac{r-2}{r-1}}\left(  \sum_{u\in V\left(  G\right)  }\frac{t^{r-1}\left(
u\right)  }{d^{2r-3}\left(  u\right)  }\right)  ^{\frac{1}{r-1}}\geq\sum_{u\in
V\left(  G\right)  }\left(  \frac{t^{r-1}\left(  u\right)  }{d^{2r-3}\left(
u\right)  }\right)  ^{\frac{1}{r-1}}d^{\frac{r-2}{r-1}}\left(  u\right)
=\sum_{u\in V\left(  G\right)  }\frac{t\left(  u\right)  }{d\left(  u\right)
},
\]
and so,%
\[
\sum_{u\in V\left(  G\right)  }\frac{t^{r-1}\left(  u\right)  }{d^{2r-3}%
\left(  u\right)  }\geq\frac{1}{\left(  2m\right)  ^{r-2}}\left(  \sum_{u\in
V\left(  G\right)  }\frac{t\left(  u\right)  }{d\left(  u\right)  }\right)
^{r-1}.
\]
Hence, from (\ref{ineq2}), we find that%
\[
n\mu_{n}\left(  G\right)  \leq-\frac{4}{\left(  r-1\right)  m^{r-2}}\left(
\sum_{u\in V\left(  G\right)  }\frac{t\left(  u\right)  }{d\left(  u\right)
}\right)  ^{r-1}.
\]

Assume the assertion of the theorem is false, that is to say
\begin{equation}
\mu_{n}\left(  G\right)  \geq-\frac{2^{r+1}m^{r}}{rn^{2r-1}},\label{oppos}%
\end{equation}
and so,
\[
\frac{4}{\left(  r-1\right)  m^{r-2}}\left(  \sum_{u\in V\left(  G\right)
}\frac{t\left(  u\right)  }{d\left(  u\right)  }\right)  ^{r-1}\leq-n\mu
_{n}\left(  G\right)  \leq\frac{2^{r+1}m^{r}}{rn^{2r-2}},
\]
implying
\[
\sum_{u\in V\left(  G\right)  }\frac{t\left(  u\right)  }{d\left(  u\right)
}\leq2\left(  \frac{r-1}{r}\right)  ^{\frac{1}{r-1}}\frac{m^{2}}{n^{2}}.
\]
Hence, in view of (\ref{ineq1}), we obtain,%
\[
\mu_{n}\left(  G\right)  \leq\left(  4\left(  \frac{r-1}{r}\right)  ^{\frac
{1}{r-1}}-4\right)  \frac{m^{2}}{n\left(  n^{2}-2m\right)  },
\]
and from (\ref{oppos}) it follows that
\[
-\frac{2^{r+1}m^{r}}{rn^{2r-1}}\leq\mu_{n}\left(  G\right)  \leq4\left(
\left(  \frac{r-1}{r}\right)  ^{\frac{1}{r-1}}-1\right)  \frac{m^{2}}{n\left(
n^{2}-2m\right)  }.
\]
Then,%
\[
\left(  1-\left(  \frac{r-1}{r}\right)  ^{\frac{1}{r-1}}\right)  \frac{m^{2}%
}{n\left(  n^{2}-2m\right)  }\leq\frac{2^{r-1}m^{r}}{rn^{2r-1}}%
\]
and so,%
\[
1-\left(  \frac{r-1}{r}\right)  ^{\frac{1}{r-1}}\leq\frac{2^{r-1}%
m^{r-2}\left(  n^{2}-2m\right)  }{rn^{2r-2}}.
\]
Since $0<m<n^{2},$ the expression $m^{r-2}n^{2}-2m^{r-1}$ attains its maximum
at $m=\frac{r-2}{2\left(  r-1\right)  }n^{2};$ thus,%
\[
1-\left(  \frac{r-1}{r}\right)  ^{\frac{1}{r-1}}<\frac{2}{r}\left(  \frac
{r-2}{r-1}\right)  ^{r-2}\frac{1}{r-1}.
\]
Hence, by Bernoulli's inequality, we see that%
\[
\frac{1}{r\left(  r-1\right)  }\leq1-\left(  \frac{r-1}{r}\right)  ^{\frac
{1}{r-1}}<\frac{2}{r}\left(  \frac{r-2}{r-1}\right)  ^{r-2}\frac{1}{r-1},
\]
a contradiction for $r\geq3.$ Therefore, assumption (\ref{oppos}) is false and
the proof is completed.
\end{proof}

We conclude with the following problem.

\begin{problem}
Find a simple lower bound on the second eigenvalue of graphs of order $n,$
size $m,$ and independence number $r.$
\end{problem}


\begin{thebibliography}{9}                                                                                                %


\bibitem {Bol98}B. Bollob\'{a}s, \emph{Modern Graph Theory,} Graduate Texts in
Mathematics, \textbf{184,} Springer-Verlag, New York (1998), xiv+394 pp.

\bibitem {BoNi04}B. Bollob\'{a}s and V. Nikiforov, Graphs and Hermitian
matrices: eigenvalue interlacing, to appear in \emph{Discr. Math.}
\end{thebibliography}
\end{document}